\documentstyle[amstex, 12pt]{amsart}

\begin{document}

\title[Complete ideals]
{Integral closedness of $MI$ and the  formula of Hoskin and Deligne
  for finitely supported complete ideals}

\author{Clare D'Cruz}
\address{Chennai Mathematical Institute, 
         G. N. Chetty Road, 
         T. Nagar, 
         Chennai 600 017 India}
\email{clare@@cmi.ac.in}
\date{\today}

\parskip = 6pt
\parindent = 10pt

\newcommand{\ncom}{\newcommand}
\ncom{\bq}{\begin{equation}}
\ncom{\eq}{\end{equation}}
\ncom{\beqn}{\begin{eqnarray*}}
\ncom{\eeqn}{\end{eqnarray*}}
\ncom{\beq}{\begin{eqnarray}}
\ncom{\eeq}{\end{eqnarray}}
\ncom{\been}{\begin{enumerate}}
\ncom{\eeen}{\end{enumerate}}
\ncom{\nno}{\nonumber}
\ncom{\hs}{\mbox{\hspace{.25cm}}}
\ncom{\rar}{\rightarrow}
\ncom{\lrar}{\longrightarrow}
\ncom{\Rar}{\Rightarrow}
\ncom{\noin}{\noindent}

\newtheorem{thm}{Theorem}[section]
\newtheorem{lemma}[thm]{Lemma}
\newtheorem{cor}[thm]{Corollary}
\newtheorem{pro}[thm]{Proposition}
\newtheorem{example}[thm]{Example}
\newtheorem{remark}[thm]{Remark}
\newtheorem{definition}[thm]{Definition}
\newtheorem{blank}[thm]{}

\ncom{\bt}{\begin{thm}}
\ncom{\et}{\end{thm}}
\ncom{\bl}{\begin{lemma}}
\ncom{\el}{\end{lemma}}
\ncom{\bco}{\begin{cor}}
\ncom{\eco}{\end{cor}}
\ncom{\bp}{\begin{pro}}
\ncom{\ep}{\end{pro}}
\ncom{\bex}{\begin{example}}
\ncom{\eex}{\end{example}}
\ncom{\brm}{\begin{remark}}
\ncom{\erm}{\end{remark}}
\ncom{\bb}{\begin{blank}}
\ncom{\eb}{\end{blank}}

\ncom{\bd}{\begin{definition}}
\ncom{\ed}{\end{definition}}

\ncom{\bc}{\begin{center}}
\ncom{\ec}{\end{center}}

\ncom{\comx}{I\!\!\!\!C}
\ncom{\zee}{$Z\!\!\!\!Z$}
\ncom{\ze}{Z\!\!\!\!Z}
\ncom{\Q}{$I\!\!\!\!Q$}
\ncom{\N}{I\!\!N}
\ncom{\sz}{\scriptsize}
\ncom{\CM}{Cohen-Macaulay }
\ncom{\sop}{system of parameters}
\ncom{\eop}{\hfill{$\Box$}}
\ncom{\tfae}{the following are equivalent:}
\ncom{\mm}{minimal multiplicity }
\ncom{\f}{\frac}
\ncom{\la}{\lambda}
\ncom{\si}{\sigma}
\ncom{\ssize}{\scriptsize}
\ncom{\al}{\alpha}
\ncom{\be}{\beta}
\ncom{\Si}{\Sigma}
\ncom{\ga}{\gamma}
\ncom{\kbar}{\overline{\kappa}}
\ncom{\bib}{\bibitem}
\ncom{\sst}{\subset}
\ncom{\sms}{\setminus}
\ncom{\seq}{\subseteq}
\ncom{\est}{\emptyset}
\ncom{\bighs}{\hspace{.5 cm}}
\ncom{\ulin}{\underline}
\ncom{\olin}{\overline}
\ncom{\bip}{\bigoplus}
\ncom{\sta}{\stackrel}
\ncom{\scl}{\succurlyeq}

\def\m{{\frak m}}
\def\spec{\mathop{\rm Spec}}

\maketitle

\section{Introduction}

It is well known that the theory of complete ideals in a two dimensional
regular local ring was founded by Zariski \cite{zar-sam}. His work
was inspired by the birational theory of linear systems on smooth
surfaces. 
Due to the complexity in higher dimension,  an higher dimensional analogue 
was not easy to obtain, as observed by Zariski and evident from the 
work of Lipman, Cutkosky, Plitant, Lejeune  etc. 
For example see \cite{lipman}, \cite{cutkosky1},  \cite{cutkosky2},
\cite{lejeune1} and \cite{cgl}.
 Their work also reveals that  a trivial
generalization of Zariski's work is not possible to obtain.

Zariski proved  that in a two dimensional regular local ring product
of  complete ideals is complete. Moreover, every complete ideal can be
uniquely factorized as a product of simple complete ideals. Both these
results do not hold true in higher dimension. A substantial amount of
work on the unique factorization has been done by Cutkosky,  Lipman and
Piltant.

Of recent interest, is the following question: If $I$ is an ideal in a
Noetherian local ring $(R,M)$, when is  $M I$ integrally closed (see
\cite{h-huneke}, \cite{corso})~? This is closely related to   the
Cohen-Macaulayness  of the  fiber cone of $I$, $F(I):=  R[It]
\otimes_R R/M = \oplus_{n \geq 0} I^n/MI^{n}$.    Hence two topics
of interest arise: the integral closedness of $MI$ and the
Cohen-Macaulayness of $F(I)$.   By an example we show that if the
dimension of the ring is at least three and $I$ is a finitely supported
complete ideal in a regular local ring,
 then the fiber cone of $I$   
need not be
Cohen-Macaulay,  in contrast to the case of dimension two where the
fiber cone of an $M$-primary complete ideal is always Cohen-Macaulay
\cite[Corollary~2.5]{c-r-v}.

Let $(R,M)$ be a regular local ring. In this paper we give necessary and sufficient conditions  in terms of
the number of generators of a finitely supported monomial complete
ideal $I$ for $MI$ to be integrally closed.  In a two dimensional
regular local ring, if $I$ is an  $M$-primary  ideal, then   $\mu(I)
\leq 1 + o(I)$, where $o(I)$ is the $M$-adic order of $I$ and $\mu(I)$
is the minimal number of generators of $I$. If $I$ is complete,  then
equality holds \cite[Lemma~3.1]{rees}. More generally, if  $I$ is an
$M$-primary complete ideal in a regular local ring $R$ of dimension $d
\geq 2$, then $\displaystyle{\mu(I) \geq  {o(I) + d-1 \choose d-1}}$
\cite{dcruz}. If $I$  is finitely supported, then,  $\mu(I) \leq
o(I)^{d-1} + d-1$ (see Theorem~\ref{generators}).  There exists
examples of finitely supported complete ideals  where the upper bound
is attained (see Example~\ref{ex-7.1}).

In this paper we prove:

\bt
\label{main-generators}
Let $I$ be a finitely supported monomial complete ideal in a regular
local ring $(R, M)$ of dimension at least three. Assume that $k = R/M$ 
is algebraically closed field. Then  $MI$ is integrally closed if 
and only if 
\beqn
\mu (I) = {o(I) + d-1 \choose d-1}.
\eeqn
\et

An essential ingredient in the proof is the formula of Hoskin and
Deligne which expresses  the  length of a finitely supported complete
ideal in terms of the order of the strict transform of $I$ (see
Theorem~\ref{main1}).  
Note that the condition
``finitely supported'' is necessary (see \cite[Lemma~(1.21.1)]{lipman}). 
This formula was known only for complete ideals of height two in 
a two dimensional regular local ring
 and   was   proved independently by Hoskin \cite{hoskin}  and Deligne
 \cite{deligne}.  
Due to several hurdles the higher dimensional analogue was not easy to
obtain. A formula in dimension three was obtained by  M. Lejeune-Jalabert
\cite{lejeune1}. 
One of the main obstacles is that the theory of complete
 ideals in a two dimensional regular local ring which was founded by
 Zariski \cite{zar-sam} does not directly  generalize  to higher 
dimension
 even though  the definitions do. In this paper we obtain the higher
 dimensional analogue of this length formula.

Let ${\olin I}$ denote the integral closure of $I$. 
A formula for  $\ell (R/{\olin{I^n}})$ for all $n \geq 1$ was obtained
by   Morales in \cite{morales2}.  In fact, his formula holds true for
a finitely generated normal $k$-algebra over an algebraically closed
field. He also gave a geometrical interpretation for the coefficients
of the Hilbert-Samuel polynomial of $I$.

For the sake of completeness we state the formula of Hoskin and Deligne:

\bt
\label{hd-2}
\cite[Theorem~5.2]{hoskin}, \cite[Theorem~2.13]{deligne}
Let $(R,M)$ be a two dimensional regular local ring with infinite 
residue field. Let $I$ be an $M$-primary complete ideal in $R$. 
  Let  ${\cal C}_I = \{Q_0, Q_1, \ldots, Q_t \}$ 
be the base points of $I$and $R_i$ the local ring at $Q_i$
$(R_0 = R)$. Then 
\beqn
  \ell \left( \f{R}{I} \right)
= \sum_{i=0}^{t} 
   {o (I^{R_i}) + 1 \choose 2}   [R_i/M_i : R/M]
\eeqn
where $M_i$ is the maximal ideal of $R_i$, $(M_0 = M)$
and $[R_i/M_i : R/M]$ denotes the degree of the field extension
$R_i/M_i \supseteq  R/M$.
\et

In higher dimension one can verify that:

\bt
\label{ineq-main}
Let $(R, M)$ be a regular local ring of dimension $d$ with infinite residue
field.  Let $I$ be a finitely supported complete ideal. 
 Let  ${\cal C}_I = \{Q_0, Q_1, \ldots, Q_t \}$ 
be the base points of $I$and $R_i$ the local ring at $Q_i$  
$(R_0 = R)$. Then  
\beqn
  \ell \left( \f{R}{ I} \right)
\leq   \sum_{i=0}^{t} 
   {o (I^{R_i}) + d-1 \choose d}   [R_i/M_i : R/M]
\eeqn
where $M_i$ is the maximal ideal of $R_i$, $(M_0 = M)$
and $[R_i/M_i : R/M]$ denotes the degree of the field extension
$R_i/M_i \supseteq  R/M$.
\et

If the dimension of the ring is at least three then the inequality in
Theorem~\ref{ineq-main} may be strict (see Example~\ref{extwo}, 
Example~\ref{ex-7.2}). The gap between the terms on the left and the
right can be estimated in terms of the length of the  right derived
functors of the direct image sheaf.

 In this paper we prove:

\bt
\label{main1}
Let $(R, M)$ be a regular local ring. Assume that $k = R/M$ is an
 algebraically closed field.
 Let $I$ be a finitely supported complete ideal.
 Let  ${\cal C}_I = \{Q_0, Q_1, \ldots, Q_t \}$ 
be the base points  of $I$ and let $\sigma: X({\cal C}_I) \lrar \spec~R$ be the
birational map obtained by a sequence of blowing up of points of
${\cal C}_I$. Then  
\beqn
   \ell \left( \f{R}{I} \right)
+ \sum_{i=1}^{d-2} (-1)^{i+1} 
  \ell ( R^i \sigma_{*} ( I {\cal O}_{X({\cal C}_I)}))
=  \sum_{i=0 }^{t} {o (I^{R_i}) + d-1 \choose d} 
\eeqn
where  $R_i$ the  local ring at the point  $Q_i$  $(R_0 = R)$ and  $M_i$
is the maximal ideal of $R_i$, $(M_0 = M)$.
\et

As a consequence of Theorem~\ref{main1} we are able to  recover the
formula for the multiplicity and the mixed-multiplicities 
of finitely supported   complete ideals  (\cite{kodiyalam}, \cite{johnston}, \cite{piltant}).

It is evident, from Theorem~\ref{ineq-main} and Theorem~\ref{main1}
that  $\sum_{i=1}^{d-2} (-1)^{i+1} \ell ( R^i
\sigma_{*} ( I {\cal O}_{X({\cal C}_I)}))  \geq 0$. This inequality
can be strict  if $I$ is not a
monomial ideal  (see Example~\ref{extwo}). 
If  $I$ is a monomial ideal, then 
$\sum_{i=1}^{d-2} (-1)^{i+1} 
\ell ( R^i \sigma_{*} ( I {\cal O}_{X({\cal C}_I)})) = 0$. 
This is a consequence of the length formula obtained by Morales
(\cite[Lemma~3]{morales3}). 
We independently obtain a closed formula for finitely  supported
monomial ideals (see Theorem~\ref{main-monomial})
which implies the same result.

We also list a few  applications. Let $(R,M)$ be a regular local ring
  of   dimension $d \geq 2$. Then the  graded ring associated
  to the   filtration  ${\cal F}  = \{ \olin {I^n} \}_{n \geq 0}$,  
$G({\cal F}) := \oplus_{n \geq 0} {\olin {I^n}}/ {\olin {I^{n+1}}}$
  is of interest. The normalization of the Rees ring
  $R[It]$ denoted by ${\olin{R[It]}}$ is the graded ring  
$\oplus_{n \geq 0}{\olin{I^n}}$. Since the depth of  $G({\cal F}) \geq
1$, $\mbox{depth }(G({\cal F})) = \mbox{depth } ({\olin{R[It]}}
  \otimes_{R}   R/I )$.

In a two dimensional regular local ring, it is well known that if $I$
is a complete ideal then  $G({\cal F}) = G(I):=\oplus_{n \geq 0}
{I^n}/ I^{n+1}$ and $G(I)$ is Cohen-Macaulay  (\cite{morales},
\cite{john-verma}).  Even in the three dimensional case $G(I)$ need
not be Cohen-Macaulay. The first example was given by Cutkosky in
\cite{dale5}. In this paper, under certain conditions we are able
to deduce the Cohen-Macaulayness of $G({\cal F})$ for a finitely
supported   complete ideal in a regular local ring of dimension three.

We end this paper with some examples which probably will put more
light on the results of this paper.

I would like to thank   J. Verma for suggesting the problem. I am also
 grateful to the   Minist{\`e}re de la recherche for financial support 
 and the Universit{\'e} de Versailles St-Quentin-en-Yvelines -   LAMA,
 Versailles, France  for local hospitality.  I am greatly indebted to
 M. Lejeune-Jalabert for her guidance and help. Special thanks are
 also due to D.~Cutkosky, V.~Kodiyalam and O.~Piltant without whom
 this paper would not be complete.  

The author also thanks the refree for his suggestions. 

\section{Preliminaries}

Let $(R, M)$ be a regular local ring. 
For all the definitions in this section we refer to 
 \cite{zar-sam}, \cite{lipman} and  \cite{cgl}.

Let $I$ be an ideal in $R$. The {\em integral closure} of $I$ denoted by
${\bar I}$ is the set 
\beqn
\{ x \in R | x^n + a_1 x^{n-1} + \cdots + a_n = 0; ~~~~~~
    a_j \in I^j, 1 \leq j \leq n \}.
\eeqn
The {\em completion} of $I$ is the ideal 
${\displaystyle
I' = \cap_{{\it v} \in S} I R_{\it v}}
$
where $S$ denotes the set of all non-trivial valuations which are
non-negative on $R$. 
Zariski proved that $I' = {\bar I}$. 
Hence the integral closure of $I$ is an ideal of $R$. 
An ideal  $I$ is {\em integrally closed or complete}
if $I = {\olin I}$.

Let $X$ be a non-singular variety and let $O$ be a point on $X$. Put
$R = {\cal O}_{X,O}$ and let $M$ be the maximal ideal of $R$. 
Let $f : X_1 \lrar  \spec ~R$ denote the blowing up of $M$.
The (first) 
{\em quadratic transforms} of $R$ are the local rings
${\cal O}_{X_1, P}$, where $P \in f^{-1}\{M \}$ is a point on $X_1$.   

A point $Q$ on a variety $Y$ is  {\it infinitely near to} 
a point $O$ on $X$,   $Q \succeq O$ in symbol,  if
\been
\item
there exists a sequence of blowing ups
$$
\sigma: Y = X_{n} \sta{f_{n}}{\lrar}
         X_{n-1}   \sta{f_{n-1}}{\lrar}
                 \cdots {\lrar}
         X_{1}    \sta{f_{1}}{\lrar}
         X_0 =    \spec ~ R;
$$
where each $f_{i+1}$, $0 \leq i \leq n-1$ is obtained by  
blowing up  $M_i$, the maximal ideal of ${\cal O}_{X_i, P_i}$,
 $P_i \in f_i^{-1}\{M_{i-1}\}\subseteq X_i$ 
($M_0 =M$, $P_0 = O$);

\item $Q \in
f_{n-1}^{-1}
\{ M_{n-1} \}$. 
\eeen

Put   $R_i := {\cal O}_{X_i, P_i}$ and $S:= R_n$. 
The sequence 
$
R= R_0 \subset R_1 
       \subset \cdots 
       \subset  R_i 
       \subset \cdots 
       \subset R_n = S
$
 is called the {\it quadratic sequence} from  $R$ to $S$.
 We   say that $S$
is {\it infinitely near} to $R$ and denote it by  $S \succeq R$.
 Given any ideal $I \subset
R$,  the {\it  transform} of $I$  in $S =R_n = {\cal O}_{X_n, P_n}$ denoted by
$I^S$ is defined
inductively as follows:
\been
\item
$I^{R_0} = I$;

\item
$I^{R_i} = M_{i-1} ^{-m_{i-1}} I^{R_{i-1}}$, where 
$m_{i-1} = o(I^{R_{i-1}})= \max \{n | I_{i-1} \seq M_{i-1}^n \}$. 
\eeen

The {\it point basis} of a non-zero ideal $I$ is a family of non-negative
integers  ${\cal B} (I) = \{ o(I^S) : S \succeq R \} $.  
We say that a point $Q \in Y$ is a {\it base point}  of $I$ if  $Q \succeq O$ 
and $o(I^S) < \infty$,  where $S = {\cal O}_{Y, Q}$.
  We say that an ideal $I$ is {\it finitely supported } if $I \not =
  0$ and if $I$  has at most  
finitely  many base points \cite{lipman}. For a finitely supported
ideal $I$ we will denote the set base points of $I$ by 
  ${\cal C}_I = \{Q_0 = O, Q_1, \ldots, Q_t\}$. 

\section{Hoskin-Deligne formula for finitely supported complete ideals}

Let $k$ be an algebraically closed field and let $X$ be a non-singular
variety of dimension at least  two. Let $O \in X$ be a 
point. Put $R ={\cal O}_{X, O}$. Let $M$ be the maximal ideal at $O$.


 The notion of  $*$-product was introduced in 
\cite{lipman}.  Let $I$ and $J$ be ideals in $R$. The $*$-product of 
 $I$ and $J$ denoted by $I * J$ is the ideal $\olin{IJ}$. An ideal $I$ is
 $*$-simple of it cannot be decomposed as a $*$-product of proper complete
 ideals  \cite{lipman}.  Notice that $I * J = {\olin{I} * {\olin {J}}}$. 

To prove Theorem~\ref{main1} we first need to consider the
first blow up. The following lemma is basically a consequence 
of Lemma~2.3 of \cite{lipman}. 
We prove it since it is the  crucial result in proving Theorem~\ref{main1}. 

\bl
\label{lemma2}
Let $f: X_1  \lrar Spec~R$ denote the blowing up of $M$. Let ${\cal I}$
be the coherent ${\cal O}_{X_1}$-ideal whose stalk at  any 
point $P \in f^{-1} \{ M \}$ is a complete ideal  and  
${\cal I}_P =  {\cal O}_{X_1,P}$ if $P \not \in f^{-1}\{M \}$. 
Then there exists a complete  ideal $I \seq R$ such that 
${\cal I}_{P} = {\olin{I^{{\cal O}_{X_1, P}}}}$, where  $P \in f^{-1} \{ M \}$.
Moreover, there exists  a  positive integer $N$ such
that for all 
$n \geq N$ 
\been  
\item
\label{two}
$
{\displaystyle 
R^j f_* (M^{ o(I) +n} {\cal I}) = 0  \hspace{.5in}
\mbox{for all } j > 0};
$\\

\item
\label{one}
$
{\displaystyle 
\ell_R \left( \f{R}{M^n * I}  \right)
=   {o(I) +n + d - 1 \choose d}
+  \sum_{P \in f^{-1}{\{ M \}}} \ell_R \left( 
   \f{{\cal O}_{X_1, P}}{{\cal I}_P}   \right)}.
$
\eeen
\el
\pf 
The existence of a complete ideal $I \subseteq R$ satisfying the
assumptions in the lemma was proved in \cite[Lemma~2.3]{lipman}. 

 From \cite[Proposition~8.5, pg. 251]{hartshorne}
\beqn
R^j f_* (M^{ o(I)+n} {\cal I}) = H^j (X_1, M^{o(I)+n} {\cal I})
\eeqn
and  for $j \geq 0$ and  all $n$ large, 
$H^j (X_1, M^{ o(I) +n} {\cal  I}) = 0$. 

We now prove (\ref{one}). 
The exact sequence

\beqn
0 \lrar {M^{o(I) + n}}{\cal I}{\cal O}_{X_1}
  \lrar  M^{o(I) + n} {\cal O}_{X_1}
  \lrar  M^n {\cal O}_{X_1} /  {\cal I}{\cal O}_{X_1}
  \lrar 0 
\eeqn

gives  the long exact sequence

\beq
\label{lemma2eq1} \nonumber
0 && \lrar H^0(X_1, M^{o(I) + n} {\cal I}{\cal O}_{X_1})
     \lrar H^0(X_1, M^{o(I) + n} {\cal O}_{X_1})
     \lrar H^0(X_1, M^n {\cal O}_{X_1} /  {\cal I}{\cal O}_{X_1}) \\
  && \lrar H^1(X_1, {M^{o(I) +n}} {\cal I}{\cal O}_{X_1})
     \lrar H^1(X_1, M^{o(I) + n}   {\cal O}_{X_1})
     \lrar \cdots.
\eeq

Now
$
 H^0(X_1, {M^{o(I) + n}} {\cal I}{\cal O}_{X_1}) 
= \olin{M^n I}
$
and 
$
  H^0(X_1, {M^{o(I) +n}}{\cal O}_{X_1}) 
=     \olin{M^{o(I) + n}} 
$ 
(see \cite[pf of Lemma~2.3]{lipman}. 
Since 
$
 H^0(X_1, M^n {\cal O}_{X_1} /  {\cal I}{\cal O}_{X_1})
=  \sum_{P \in f^{-1} \{M\}} 
     \ell \left({\cal O}_{X_1, P}/{\cal I}_{P} \right)
$
and for all $n \gg 0$, 
$
  H^1(X_1, {M^{o(I) + n} {\cal I}}{\cal O}_{X_1}) = 0
$ 
(see  \cite[Theorem~5.2, pg 228]{hartshorne}),
plugging these in  the exact sequence (\ref{lemma2eq1}) gives the result. 
\eop

\bp
\cite[Proposition~2.1, Corollary~2.2]{lipman}
Let $S \succeq R$. Then there exists a unique $*$-simple complete ideal
$p_{RS} \subseteq R$ satisfying the following properties:
\been
\item
$R/ p_{RS}$ is artinian;

\item
$p_{RR}$ is the maximal ideal of $R$;

\item
For all regular local rings $T$ with $S \succeq T \succeq R$, 
$p_{TS} = \olin{(p_{RS})^T}$.
\eeen 
\ep

We now apply Lemma~\ref{lemma2} recursively to a sequence of point
blowing ups.

\bl 
\label{lemma3}
Let $I$ be a finitely supported ideal and  ${\cal C}_I= \{Q_0 = O,
Q_1, \ldots, Q_t\}$  the base points  of  $I$. 
Let  $R_i$ denote the local ring at $Q_i \in {\cal C}_I$. 
Put  $m_i = o({I^{R_i}})$. There exists integers $N_0, 
  \ldots, N_t$ such that for all $n_i \geq N_i$, 
\been

\item
\label{five}
$
{\displaystyle
R^{j} f_* ({{{\prod}^*}_{i=0}^{t}}~(p_{RR_i}^{n_i}* {I}){\cal O}_{X_1({\cal C}_I)}) = 0}$
for all $j \geq 0$;
\\

\item
\label{four}
$
{\displaystyle
  \ell \left( 
  \f{R}{{{\prod^*}_{i=0}^{t}}~{p_{RR_i}^{n_i}} *I}  \right)
= \sum_{i=0}^{t} {m_i + \sum_{R_j \succeq R_i} n_j + d-1 \choose d}}.
$
\eeen
Here  ${{\prod}^*}$ denotes $*$-product.
\el
\pf The first part follows from  \cite[Proposition~8.5,
pg. 251]{hartshorne}. 

We prove (\ref{four}) by induction on $t = |{\cal C}_I|-1$. If  $t = 0$,
then ${\olin I} = M^n$ for some $n >0$ and (\ref{four})
follows trivially. 

Let $t >0$. Let  $f : X_1 \lrar Spec~R$ denote the blowing up of
$Spec~R$ at $M$. Let ${\cal I}$ be the ${\cal O}_{X_1}$-ideal sheaf whose
stalk at every point $Q_i \in {\cal C}_I \cap X_1$ is  
\beqn
{\cal I}_{Q_i} = {\olin {I^{R_i}}} =
\olin{M^{-o(I)} I  {\cal O}_{X_1, Q_i}}
\eeqn
 and   ${\cal I}_{Q} =
{\cal O}_{X_1, Q}$ for
 $Q \in X_1 \sms (X_1 \cap {\cal C}_{{\cal I}})$.

Let ${\cal C}_{{\cal  I}_{Q_i}} = \{Q_{i(0)}= Q_i, \ldots, Q_{i(s_i)} \}$. 
For all $Q_i \in X_1 \cap {\cal C}_I$,  ${\cal C}_{{\cal I}_{Q_i}} 
\subset
{\cal C}_{\cal I}$   and hence  
 $|{\cal C}_{{\cal  I}_{Q_i}}| < |{\cal C}_{\cal I}|$.

By induction 
hypothesis, for each $Q_i \in  {\cal C}_I \cap X_1$, 
there exists integers $N_{i(0)}, \ldots, N_{i(s_i)}$ 
such that for all $n_{i(j)} \geq N_{i(j)}$,  $0 \leq j \leq s_i$,  
\beq
\label{lemma3eq1}
  \ell \left( \f{{\cal O}_{X_1, Q_i}}
                {\displaystyle {{\prod}^*}_{j=0}^{s_i}~  
                 p_{R_i R_{i(j)}}^{n_{i(j)}} * {\cal I}_{Q_i}} \right)
= \sum_{j=0}^{s_i}
  {m_{i(j)} + \sum_{R_k \succeq R_{i(j)}} n_k + d-1 \choose d}.
\eeq

Fix $n_{i} \geq N_{i}$, $0 \leq i \leq t$. 
Let ${\cal J}(n_1, \ldots, n_t; I)$ 
be the ${\cal O}_{X_1}$- ideal sheaf whose stalk at
every point $Q_i \in {\cal C}_{I} \cap X_1$ is 
\beqn
   {\cal J}( n_1, \ldots, n_t; I)_{Q_i}
:=  \left( {{\prod}^*}_{j=0}^{s_i}~
    p_{R_i R_{i(j)}}^{n_{i(j)}} * {\cal I}_{Q_i} \right)
\eeqn
and 
$
{{\cal J}}(n_1, \ldots, n_{t}; I)_{Q} = {\cal O}_{X_1, Q}
$
 for   $Q \in X_1 \sms (X_1 \cap {\cal C}_{{\cal I}})$.
Let $n_0 \geq 0$ and for each $n_0$ let 
$$
J(n_0, n_1, \ldots, n_t; I) 
=  p_{R R}^{n_0} *  {{\prod}^*}_{j=1}^{t} p_{R R_j}^{n_j} * { I}.$$ 
Here $p_{RR} =
M$, the maximal ideal of $R$.  
Then
$$
{\olin{ M^{-o(J(n_0, n_1, \ldots, n_t; I))}
             J(n_0, n_1, \ldots, n_t;I){ {\cal O}_{X_1, Q_i}}}}
= {{\cal J}}(n_1, \ldots, n_{t}; I)_{Q_i} 
$$
 at all points $Q_i \in X_1 \cap {\cal C}_I$.  

Now by  Lemma~\ref{lemma2} 
there exists and integer $N_0$ such that for each $n_0 \geq N_0$, 
\beq
\label{lemma3eq2} \nno
&&   \ell_R \left( 
     \f{R}
       {J(n_0, n_1 \ldots, n_t; I)} \right) \\ \nno
&=&  \ell_R \left( \f{R}{p_{RR}^{ o(J (n_0, n_1, \ldots, n_t;I))}} \right)
+    \sum_{Q_i \in {\cal C}_I \cap X_1}
     \ell_{{\cal O}_{X_1, Q_i}} \left(  
     \f{{\cal O}_{X_1, Q_i}}
       {{{\cal J}}(n_1, \ldots, n_{t};I)_{Q_i}} \right) \\
&=&    {o(I) + n_0 + \cdots + n_t +   d-1 \choose d}
+    \sum_{Q_i \in {\cal C}_I \cap X_1}
     \ell_{{\cal O}_{X_1, Q_i}} \left(  
     \f{{\cal O}_{X_1, Q_i}}
       {{{\cal J}}(n_1, \ldots, n_{t}; I)_{Q_i}} \right). 
\eeq
 Substituting  
(\ref{lemma3eq1}) in (\ref{lemma3eq2}) proves the lemma. 
\eop

\vspace{.2in}

Let ${\cal C}_I = \{Q_0 = O, Q_1, \ldots, Q_t\}$ denote the base
points of a  finitely supported ideal $I$.  Let $X({\cal C}_I)$
denote the variety obtained by blowing up $X_t$ at  $Q_t$.  Let $E_i$
be the exceptional divisor obtained by blowing up $Q_i$ and let
$E_{i}^*$ denote the exceptional divisor in $X_h, i \leq h \leq n+1$.  
Let ${\cal A}_I = \{ {\cal C}_I, B(I) \}$ where 
${\cal B}(I) = \{ m_0, \ldots, m_t \}$ is the  point basis of $I$. 
Then ${\cal D}({\cal A}_I)
= \sum_{i=0}^{t}{m_i E_{i}^*}$ is the divisor associated to the ideal sheaf
$I {\cal O}_{X({\cal C}_I)}$. 
Let $R_i$ be the regular local ring at $Q_i$. 
Then the exceptional divisor corresponding 
to 
$(\prod_{i=0}^{t} p_{RR_i}^{n_i} *  I) 
{\cal O}_{X({\cal C}_I)}$ is $\sum_{i=0}^{t}{h_i E_{i}^*}$ where
$h_i  =  m_i  + \sum_{R_j \succeq R_i} n_j$. 

\vspace{.5in}

\noin
{\bf Proof of Theorem~\ref{main1}:}
Denote ${\cal L}_i = {\cal O}_{X({\cal C}_I)} (-E^*_i) $
 and let  $h_i  =  m_i  + \sum_{R_j \succeq R_i} n_j$. 
We have
\beq
\label{six} \nno
&&  \chi \left( {\cal L}_0^{\otimes h_0} \otimes \cdots  
                {\cal L}_t^{\otimes h_t} \right)  \\ \nno
&=&   \ell \left( 
    \f{R}
      {\prod_{i=0}^{t} p_{RR_i}^{n_i} *  I} \right) 
+   \sum_{i=1}^{d-1} (-1)^{i+1} 
    \ell \left(R^i {\sigma}_* 
    \left(  {\prod_{i=0}^{t} \left( p_{RR_i}^{n_i}* I \right) } 
               {\cal O}_{X({\cal C}_I)} \right) \right)  \\ \nno
&&  \hspace{3in}   (\mbox{by \cite[Theorem~1.4]{morales}})  \\ \nno
&=& \sum_{i=0}^t {m_i + \sum_{R_j \succeq R_i} n_j + d-1 \choose d}
+   \sum_{i=1}^{d-1} (-1)^{i+1} 
    \ell \left( R^i {\sigma}_* 
    \left(  {\prod_{i=0}^{t} \left( p_{RR_i}^{n_i}* I \right)  } \right) 
{\cal O}_{X({\cal C}_I)} \right)   \\
&& \hspace{3in} (\mbox{ by  Lemma~\ref{lemma3}})
\eeq
for all $n_0, \ldots, n_t \gg 0$. 

On the other hand, for all non-negative 
integers $r_0, \ldots, r_t$ there exists
rational numbers  $a_{i_0 \ldots i_t}$ such that 
([Theorem~9.1]\cite{snapper})  
\beq
\label{eight}
   \chi \left( {\cal L}_0^{\otimes r_0} \otimes \cdots  
                {\cal L}_t^{\otimes r_t} \right)
= \sum_{i_0 + \cdots + i_t \leq d}
   a_{i_0 \ldots i_t}
  {r_0 + i_0 \choose i_0} \cdots  {r_t + i_t \choose i_t}
\eeq

If we put $r_i = h_i$ in (\ref{eight}), then for $n_0, \ldots, n_t \gg 0$  
 large the polynomials in 
(\ref{six}) and (\ref{eight}) agree. This gives
\beqn  
a_{i_0 \ldots i_t}  = \left\{ \begin{array}{ll}
                   1 & {\em if } \;
                  (i_0, \ldots, i_t) =  (0, \ldots,d, \ldots, 0)\\
                   -1 & {\em if } \;
                  (i_0, \ldots, i_t) =  (0, \ldots,d-1, \ldots, 0)\\ 
                   0 & \mbox{otherwise}          \end{array}
                                           \right. .
\eeqn                    
Hence
\beq
\label{twelve} \nno
 \chi \left( {\cal L}_0^{ \otimes h_0} 
     \otimes \cdots
      \otimes {\cal L}_s^{ \otimes h_t} 
                                            \right) 
= \sum_{i=0}^t \left[
  {h_i + d  \choose d} -  {h_i + d -1 \choose d-1}
\right]
\eeq
for all values of $h_0, \ldots, h_t \geq 0$. Hence (\ref{eight}) is true
     for all values of $n_0, n_1, \ldots, n_t$. 
If we put $n_0 = n_1 = \cdots = n_t = 0$ in (\ref{six}) and (\ref{eight}) we get
\beqn
  \ell \left( \f{R}{I} \right)
=   \sum_{i=0}^t {m_i + d-1 \choose d}
-  \sum_{i=1}^{d-1} (-1)^{i+1} R^i {\sigma}_*( I{\cal O}_{X({\cal C}_I)}) .
\eeqn 
It remains to show that 
$ R^{d-1} {\sigma}_* (  I  {\cal O}_{X({\cal C}_I)}) =0 $. 
Consider the exact sequence
\beqn
      0 
\lrar {\cal F} 
\lrar {\cal O}_{X({\cal C}_I)} 
\lrar I {\cal O}_{X({\cal C}_I)}
\lrar 0
\eeqn
where ${\cal F}$ is a coherent ${\cal O}_{X({\cal C}_I)}$-module. This  
     gives the exact sequence
\beqn
      R^{d-1}\sigma_* \left( {\cal O}_{X({\cal C}_I)} \right)
\lrar R^{d-1}\sigma_* \left( I {\cal O}_{X({\cal C}_I)}\right)
\lrar R^{d}\sigma_* {\cal F}
\lrar R^{d}\sigma_* \left( {\cal O}_{X({\cal C}_I)} \right).
\eeqn
Since $\sigma$ is the  composition of  sequence of blowing ups and $X$
     is non-singular, 
$R^{i}\sigma_* \left( {\cal O}_{X({\cal C}_I)} \right)= 0$
for all $i \geq 1$. And   
$ R^{d} {\sigma}_* {\cal F} = 0$ since
     $\sigma^{-1}\{M \}$ has dimension $d-1$. Hence, 
$R^{d-1}\sigma_* \left( I {\cal O}_{X({\cal C}_I)} \right)=0$. 
\eop

\section{The length formula for finitely supported monomial ideals}

The following result was proved by Morales:

\bt
\cite[Lemma~6]{morales3}
\label{length}
Let $I_1, \ldots, I_t$ be  monomial
ideals of height $d$  in $R = k[x_1, \ldots, x_d]$. Then  for all
nonnegative integers 
$r_1, \ldots, r_t$,  $\ell (R/ \olin{I_1^{r_1} \cdots I_t^{r_t}})$ is a
polynomial of degree $d$ in $r_1, \ldots, r_t$. 
\et

\brm
Let $d \geq 2$ and  let $R = k[x_1, \ldots, x_d]$ be a polynomial
ring in $d$ variables over a field $k$ and let $I \seq R$ be a
complete ideal of  height $d$. 
\been
\item
Since  $M  = (x_1, \ldots, x_d)$ is the only maximal ideal which
contains $I$, $\ell (R/I) = \ell (R_M/I_M)$. 

\item
The first quadratic transform of $R_M$ are  the  local rings 
$S_i = R[M/x_i]_{M/x_i}$. Each $T_i = R[M/x_i]$ is a polynomial ring
in $d$ variables over the field $k$. 

\item
Let $I$ be a monomial ideal in $R$ and $o(I) = \max \{n | I \seq M^n \}$. 
Then $I^{T_i} := x_{i}^{-o(I)} I$
is a monomial ideal and $I^{S_i} = I^{T_i}_{M_i}$, where $M_i$ is the
maximal ideal of $S_i$. 

\item
If $I$ is a finitely supported ideal, then $\ell (T_i/I^{T_i})$ is
finite and  $\ell (T_i/I^{T_i}) = \ell (S_i/I^{S_i})$. 
\eeen
\erm

Hence it is possible to use the theory of length of ideals in  local rings to
obtain our main result.

\bt
\label{main-monomial}
Let $I$ be a finitely supported monomial ideal in a polynomial ring
$R = k[x_1, \ldots, x_d]$. Let $M = (x_1, \ldots, x_d)$.   Then
\beqn
   \ell \left( R/ {\olin I} \right)
=  \sum_{S  \succeq R} {o (I^S) + d-1 \choose d} [S/M_S : R/M],
\eeqn
where $[S/M_S : R/M]$ denotes the degree of the field extension
$S/M_S \supseteq  R/M$. Here $M_S$ is the maximal ideal of $S$
and $M_R = M$.
\et

\pf 
We prove the theorem by induction on
$\ell_R(R/{\olin I})$. If $\ell (R/{\olin I}) = 1$, then 
${\olin I}=M$ and the result is trivially true. 

Let $\ell (R/{\olin I}) > 1$. Let  $S_i$ denote  the first quadratic transform
of $R$ and $M_i$ the maximal ideal of $S_i$. By Theorem~\ref{length}, 
$ \ell (S_i / \olin{(I^{S_i})^s})$ is a polynomial in $s$ for all $s
\geq 0$. Fix $s >0$. Then by Lemma~\ref{lemma2}, there exists an
integer $r_s$ so that for all $r \geq r_s$, 
\beq
\label{induction} \nno
    \ell \left( 
       R/{\olin{M^r I^s}} \right)
&=&  {r + s~o(I) + d-1 \choose d} 
+   \sum_{S_i} \ell_R \left(  S_i /{\olin{(I^{S_i})^s}} \right)  \\ 
&=&  {r + s~o(I) + d-1 \choose d} 
+   \sum_{S_i} \ell_{S_i} \left( 
      S_i /{\olin{(I^{S_i})^s}} \right) [S_i/M_i : R/M].
\eeq

But  
$\ell \left( R/ \olin{M^r I^s} \right)$,  is a polynomial in $r,s$ for
all $r,s \geq 0$. 
Put $r=0$ and $s=1$ in (\ref{induction}). The we have
\beqn
  \ell \left( R/ {\olin{I}} \right)
=  { o(I) + d-1 \choose d} 
+ \sum_{S_i} \ell_{S_i} \left( S_i / {\olin{I^{S_i}}} \right) [S_i/M_i : R/M].
\eeqn
Since 
$\ell_{S_i} (S_i/ {\olin {I^{S_i}}}) < \ell_R (R/ {\olin I})$, 
the result follows by induction hypothesis. 
\qed

\section{Mixed-multiplicities and the integral closedness of $MI$}

Let $(R,M)$ be a normal local ring of dimension $d$ with infinite 
residue field. It is well known that if $I_1, \ldots, I_g$ are $M$-primary
 ideals in   $R$, then for all $n_1, \ldots, n_g \gg 0$, 

 \beqn
    \ell \left( R/{ {\olin{I_1^{n_1} \cdots I_g^{n_g}}}}\right)
&=& \sum_{i_1+ \cdots + i_g \leq d}
    e_{i_1, \ldots, i_g}(I_1, \ldots, I_g)
    {n_1  + i_1  \choose i_1} \cdots 
    {n_1  + i_g  \choose i_g}
\eeqn
where $    e_{i_1, \ldots, i_g}(I_1, \ldots, I_g)$ are integers. 
For $i_1+ \cdots + i_g = d$,  $e_{i_1 \ldots i_g}$ are the mixed
 multiplicities of the ideals $I_1, \ldots, I_g$ (see \cite{teissier}, 
\cite{rees}) 
We let  ${\ulin{i}}$ denote the multi-index  $\{i_1, \ldots, i_g \}$, 
such that $\sum i_j = d$.

For the rest of this section we
will assume that $k$ is an  algebraically closed field. 
There is   a more precise formula for the mixed multiplicities of
finitely supported monomial ideals.

\bt
\label{mixed-mult1}
\cite[Corollary~3.14]{johnston}, \cite[Proposition~2.2]{piltant}
Let $I_1, \ldots, I_g$ be finitely supported  ideals in a
regular local ring $(R,M)$ of dimension $d$. Then 
\beqn
   e_{\ulin{i}}(I_1, \cdots, I_g)
= \sum_{S\succeq R} (o(I_1^{S}))^{i_1} \cdots   
                    (o(I_g^{S}))^{i_g}
\eeqn

\et
\pf Imitating the proof of Lemma~\ref{lemma2} for several ideals, we
get that for fixed   $n_1, \ldots, n_g$ there exists an  $N_0
   \geq 0$ depending on $n_1, \ldots, n_g$ such that  for all
$n_0 \geq N_0$
 \beqn
&&    \ell \left( 
    R/
      {\olin{M^{n_0}I_1^{n_1} \cdots I_g^{n_g}}} \right)\\
&=&  {n_0 + n_1o(I_1) + \cdots + n_go(I_g) + d - 1 \choose d}
+  \sum_{S_i} \ell_{S_i} \left( 
   S/{\olin{(I_1^{S_i})^{n_1} \cdots (I_g^{S_i})^{n_g} }}\right).
\eeqn
where $S_i$ are the first quadratic transform of $R$. 
Choose $n_0, n_1, \ldots, n_g \gg 0$ so that both 
$
    \ell \left( 
    R/
      {\olin{M^{n_0}I_1^{n_1} \cdots I_g^{n_g}}} \right)
$
and 
$ 
   \ell_{S_i} \left(
     S/{\olin{(I_1^{S_i})^{n_1} \cdots (I_g^{S_i})^{n_g}}} \right)
$
are polynomials. 
Now use the fact that 
$  e_{{0, i_1, \cdots, i_g}}(M, I_1, \cdots, I_g)
=  e_{{i_1, \cdots, i_g}}(I_1, \cdots, I_g)$. 
\eop

It is  well known that $e_{0,\ldots,d,\ldots,0}(I_1, \cdots, I_g) 
= e(I_i)$ where $(0,\ldots,d,\ldots,0)$ denotes the tuple where 
$d$ is at the i-th spot 
(\cite{rees1}). 
When we deal with two ideals, we will use the notation
$e_{i}(I|J) := e_{d-i, i }(IJ)$.  Note that 
$e_{i}(I|J)= e_{i}(\olin{I}|\olin{J})$. 

In a  regular local ring of dimension at least two, for any $M$-primary
complete ideal we have $e_1(M|I) = o(I)$ \cite[Theorem~4.1]{verma1},
\cite[Lemma~1.1]{verma2}.  We have an analogue of this result  for  finitely
supported complete ideals in regular local rings of 
dimension at least two.

As an immediate consequence we have:

\bco
\label{mixed-mult}
Let $I$ be finitely supported ideal in a
regular local ring $(R,M)$ of dimension $d \geq 2$. 
Then  
\been
\item
${\displaystyle e(I) = \sum_{S\succeq R} o(I^{S})},$\\

\item
$ {\displaystyle e_i(M | I) = o(I)^i \mbox{ for } 1 \leq i \leq d-1. }$
\eeen
\eco
\pf
Both (1)  and (2) follow directly from Theorem~\ref{mixed-mult1}. 
\eop

It is well known that for every $M$-primary complete ideal in a two  
dimensional regular local ring, $\mu(I) = 1 + o(I)$. We have a
generalization of this result.

\bt
\label{generators}
Let $I$ be a finitely supported  ideal in a regular local ring 
$(R,M)$ of dimension $d \geq 2$. Assume that $R/M$ is an algebraically closed
field. Then 
\been
\item
\label{gen-one}
$
{\displaystyle
\ell \left(\olin {I}/ \olin {M I} \right) 
\geq {o(I) + d-1 \choose d-1}}
$
and equality holds if and only if  
$$
\sum_{i=1}^{d-2} (-1)^{i+1} 
  \ell ( R^i \sigma_{*} ( I {\cal O}_{X({\cal C}_I)}))=0.
$$ 
 \\

\item
Assume that   $I$  integrally closed.
\been
\item
$MI$ is integrally closed of and only if 
\beqn
  \mu(I) 
=  {o(I) + d-1 \choose d-1}
+ \sum_{i=1}^{d-2} (-1)^{i+1} \left[
  \ell ( R^i \sigma_{*} (I {\cal O}_{X({\cal C}_I)}))
 - 
   \ell ( R^i \sigma_{*} (M I {\cal O}_{X({\cal C}_I)})) \right].
\eeqn
\\

\item
$
{\displaystyle
  \mu(I) 
\geq  {o(I) + d-1 \choose d-1}}
$
and equality holds if and only if $MI$ is integrally closed and 
$ 
{\displaystyle
\sum_{i=1}^{d-2} (-1)^{i+1} 
  \ell ( R^i \sigma_{*} ( I {\cal O}_{X({\cal C}_I)}))} = 0 
$
\\
\eeen

\item
\label{generators-3}
For $i=1, \cdots, d-1$,
$$
     {o(I) + d-1 \choose d-1} 
\leq \mu(I) 
\leq  o(I)^{d-i} + d-i + (i-1)~ \ell (R/I) .  
$$

In particular when $i=1$, 
$
\mu(I)  \leq o(I)^{d-1} + d-1. 
$
\eeen
\et 
\pf
We can choose an element $x \in M \sms M^2$ such that $\olin{MI} :
xR = I$ \cite[Lemma~3.1]{dcruz}. 
The exact sequence 
\beqn
0 \lrar            \f{\olin I}  {\olin{MI}}
  \lrar            \f{M^{o(I)}}      {\olin{MI}}
  \sta{.x}{\lrar}  \f{M^{o(I)}}      {\olin{MI}}
  \lrar            \f{M^{o(I)}}      {\olin{MI} + xM^{o(I)}}
  \lrar 0
\eeqn
gives
\beqn
      \ell \left( \f{\olin I}  {\olin{MI}} \right)
=     \ell \left( \f{M^{o(I)}}{\olin{MI} + xM^{o(I)}} \right)
\geq  \ell \left( \f{M^{o(I)}}{M^{o(I)+1}} \right)
=                        {o(I) + d-1 \choose d-1}. 
\eeqn
This proves (\ref{gen-one}). 

 From Theorem~\ref{main1} it follows that 
\beqn
    \ell (R/ \olin {MI}) 
&=& \ell (R/ \olin {I})
-    {o(I) + d-1 \choose d}
+    {o(I) + d \choose d-1}\\
&&
+   \sum_{i=1}^{d-2} (-1)^{i+1} 
    \ell ( R^i \sigma_{*} ( I {\cal O}_{X({\cal C}_I)}))
-   \sum_{i=1}^{d-2} (-1)^{i+1} 
    \ell ( R^i \sigma_{*} ( MI {\cal O}_{X({\cal C}_I)}))\\
&=& \ell (R/ \olin {I})+ {o(I) + d-1 \choose d-1}\\
&&+ \sum_{i=1}^{d-2} (-1)^{i+1} 
    \ell ( R^i \sigma_{*} ( I {\cal O}_{X({\cal C}_I)}))
-   \sum_{i=1}^{d-2} (-1)^{i+1} 
    \ell ( R^i \sigma_{*} ( MI {\cal O}_{X({\cal C}_I)})). 
\eeqn

Hence

\beqn
    \ell \left( \f{I}{\olin {MI}} \right)
&=& {o(I) + d-1 \choose d-1}\\
&&+ \sum_{i=1}^{d-2} (-1)^{i+1} 
    \ell ( R^i \sigma_{*} ( I {\cal O}_{X({\cal C}_I)}))
-   \sum_{i=1}^{d-2} (-1)^{i+1} 
    \ell ( R^i \sigma_{*} ( MI {\cal O}_{X({\cal C}_I)})). 
\eeqn

Applying (\ref{gen-one}) we get 
$$
\sum_{i=1}^{d-2} (-1)^{i+1} 
    \ell ( R^i \sigma_{*} ( I {\cal O}_{X({\cal C}_I)}))
\geq   \sum_{i=1}^{d-2} (-1)^{i+1} 
    \ell ( R^i \sigma_{*} ( MI {\cal O}_{X({\cal C}_I)})).
$$
Recursively we can prove:
\beqn
       \sum_{i=1}^{d-2} (-1)^{i+1} 
       \ell ( R^i \sigma_{*} ( I {\cal O}_{X({\cal C}_I)}))
&& \geq   \sum_{i=1}^{d-2} (-1)^{i+1} 
       \ell ( R^i \sigma_{*} ( MI {\cal O}_{X({\cal C}_I)})) \\
&& \geq \\
&& \vdots \\
&& \geq \\
&& \vdots \\
&& \geq
       \sum_{i=1}^{d-2} (-1)^{i+1} 
       \ell ( R^i \sigma_{*} ( M^nI {\cal O}_{X({\cal C}_I)})).
\eeqn
But 
\beqn
 \sum_{i=1}^{d-2} (-1)^{i+1} 
    \ell ( R^i \sigma_{*} ( M^nI {\cal O}_{X({\cal C}_I)})) = 0
\eeqn
for $n \gg 0$. Now, equality holds if and only if 
\beqn
 \sum_{i=1}^{d-2} (-1)^{i+1} 
    \ell ( R^i \sigma_{*} ( M^nI {\cal O}_{X({\cal C}_I)})) = 0
\eeqn
for all $n \geq 0$.

If $MI$ is integrally closed if and only if  $MI =
\olin{MI}$. Now apply (1). This proves 2(a).
2(b) follows from that fact  that $\mu(I) \geq \ell(I/ {\olin{MI}})$. 

From \cite[Theorem~2.2]{clare-verma} it follows that, for all 
$i=1, \ldots, d-1$, 
\beqn
\mu(I)
&\leq& d-i + (i-1) \ell \left( \f{R}{I} \right)
+      e_{d-i} (M |I)\\
&=&    d-i + (i-1) \ell \left( \f{R}{I} \right)
+      o(I)^i 
\hspace{.5in} \mbox{[by Corollary~\ref{mixed-mult}]}.
\eeqn
This proves(\ref{generators-3}). 
\qed

As an immediate consequence we have:

\bt
\label{generators1}
If in addition to the conditions in Theorem~\ref{generators},  
$I$ is a
 monomial   ideal and $M$ is the maximal homogeneous ideal in 
$k[x_1, \cdots, x_d]$,  then 
\been
\item
$
{\displaystyle
   \ell \left(\olin {I}/ \olin {M I} \right) 
=   {o(I) + d-1 \choose d-1}}.
$

\item
Let $d \geq 3$.  If $I$ is integrally closed, then 
$MI$ is integrally closed of and only if 
\beqn
  \mu(I) 
=  {o(I) + d-1 \choose d-1}.
\eeqn
\eeen
\et
\pf If $I$ is a monomial ideal then $MI$ is also a monomial ideal. Now,
 for any monomial ideal  $I$, $\sum_{i=1}^{d-2} (-1)^{i+1} 
    \ell ( R^i \sigma_{*} (I {\cal O}_{X({\cal C}_I)})) = 0$.
\qed

\brm
{\em
If  $I$ is a complete $M$-primary ideal in a ring of dimension at least
three and if $I$ is not finitely supported, then  
Corollary~\ref{mixed-mult} and Theorem~\ref{generators}(3) 
may not hold true.  }
\erm

\section{The associated graded ring and the Rees ring}

Let $(R,M)$ be a Noetherian local ring of positive dimension $d$. Let  $I$
be  an $M$-primary ideal in $R$.  
Here $R(I)$ and $G(I)$ will denote the ordinary Rees
ring and the associated graded ring respectively.
The filtration ${\cal F} =
\{\olin {I^n} \}_{n \geq 0}$ is a Hilbert filtration. The  Rees
ring of ${\cal   F}$,  $R({\cal F}) := \oplus \olin {I^n}$ (resp. the
associated graded ring of $G({\cal F}) := \oplus \olin {I^n} /  \olin
{I^{n+1}}$),  is a   graded ring which is  Noetherian and 
$R({\cal F})$ (resp. $G({\cal F})$) is a finite $R(I)$ (resp. $G(I)$) module. 

An ideal $J \seq \olin I$ is a reduction of ${\cal F}$
if $J \olin    {I^n} = \olin {I^{n+1}}$ for all $n \gg 0$
\cite{north-rees}.  A minimal reduction of ${\cal F}$ is a  reduction
of ${\cal F}$ which is  minimal with respect to containment.   

Since $R({\cal F})$ is a finite $R(I)$-module,  any minimal reduction
of $I$ is also a minimal reduction of ${\cal F}$. By
\cite{north-rees}, minimal reductions always exist and if the residue
field $R/M$ is infinite, then any minimal reduction of $I$ is generated
by $d$ elements.
For any minimal reduction $J$ of ${\cal F}$ we we set $r_J({\cal F}) =
sup \{ n \in \ze | J  \olin{I^{n-1}} \not = \olin{I^n} \}$

The reduction number of ${\cal F}$, denoted by $r({\cal F})$ is
defined to be the least $r_J ({\cal F})$ over all possible minimal
reductions of $J$ of ${\cal F}$. 
For any $M$-primary ideal $I$ in a local ring $(R,M)$, let $J =
 (x_1, \cdots, x_d)$ be a minimal reduction of $I$ and let  
 $C.(n):= C.(J, {\cal F}, n)$
denote the complex  
\beqn
       0 
\lrar \f{R}{\olin{I^{n-d}}} 
\lrar \cdots  
\lrar \f{R}{\olin{I^{n}}}
\lrar 0
\eeqn
where the maps are those of the Koszul complex of $R$ with respect to
$x_1, \ldots, x_d$. For details see \cite{marley}. 
Let $H_i (C.( n))$ denote the $i$-th-homology.
Let 
\beqn
h_i (J,{\cal F} ) 
= \sum_{n \geq 1} \ell (H_i (C.( n)))
= \sum_{n \geq i+1} \ell (H_{i} (C.( n))) 
\eeqn
since $H_i (C.( n)) =0$ for $n \leq i$. 

If $I$ is an $M$-primary ideal in a regular local ring of dimension two, then
$r(I) \leq 1$ \cite{huneke} 
and hence $G(I)$ is Cohen-Macaulay \cite{john-verma} and hence 
the Rees ring is Cohen-Macaulay by  \cite{goto-shimoda}.

In higher dimension, if $I$ is not a finitely supported complete ideal then it
is easy to see that both the Rees ring  $R({\cal F})$ and the
associated graded ring  $G({\cal F})$ need not be  Cohen-Macaulay.

\bl
\label{dim3}
Let $I$ be a finitely supported 
  $M$-primary  ideal in a regular local ring $(R, M)$
of dimension at least three. Assume that for all $n \geq 1$
\beqn
   \ell \left( \f{R}{\olin{I^n}} \right)
=  \sum_{S  \succeq R} {n~o(I^S) + 2 \choose 3} [S/M_S : R/M] 
\eeqn
where $M_S$ is the maximal ideal of $S$, $M_R = M$. 
Then 
\been
\item
$depth~G({\cal F}) \geq 2$ if and only if  $r({\cal F}) \leq 2$. 

\item
$G({\cal F})$ is \CM if and only if 
\beqn
 \ell \left( \f{ J +  \olin{I^2}}{J} \right)
=  \sum_{S \succeq R}{o(I_S) \choose 3} [S/M_S : R/M]
\eeqn
where $J$ is a minimal reduction of $I$. 

\item
If $o(I)  \leq 2$, then $G({\cal F})$ is \CM.

\item
If $o(I) \geq 3$, the $r_J({\cal F}) \geq 2$ for any
minimal reduction $J$ of $I$. 
\eeen
\el
\pf 
First note that for all $n \geq 3$, 
\beqn
  \ell \left( \f{R}{\olin{I^{n}}} \right)
-   3 \ell \left( \f{R}{\olin{I^{n-1}}} \right)
+   3 \ell \left( \f{R}{\olin{I^{n-2}}} \right)
-     \ell \left( \f{R}{\olin{I^{n-3}}} \right)
= e(I).
\eeqn

Since  $\sum_{i \geq 2}h_i (J,{\cal F})
\geq 0$  \cite[Theorem~3.7]{marley}, 
\beqn
&&  \sum_{i \geq 2}h_i (J,{\cal F})\\
&=&  h_2 (J,{\cal F} ) - h_3 (J,{\cal F} )\\
&=& \sum_{n \geq 3}  
    \ell \left( H_2 (C.( n)) \right) 
-   \ell \left( H_3 (C.( n)) \right) \\
&=& \sum_{n \geq 3}   \left[
    \ell \left( H_1 (C.(n)) \right)
-   \ell \left( H_0 (C.(n)) \right)
+     \ell \left( \f{R}{\olin{I^{n}}} \right)
-   3 \ell \left( \f{R}{\olin{I^{n-1}}} \right)
+   3 \ell \left( \f{R}{\olin{I^{n-2}}} \right)
-     \ell \left( \f{R}{\olin{I^{n-3}}} \right) \right]\\
&=& \sum_{n \geq 3}   \left[
    \ell \left( \
    \f{J \cap {\olin{I^n}}}
              {J \olin{I^{n-1}}}  \right)
-   \ell \left( \f{R}{J + \olin{I^n}}\right)
+    e(J) \right] \\
&=& \sum_{n \geq 3} 
    \ell \left( \f{\olin{I^n}}{J {\olin{I^{n-1}}}} \right)\\
&\geq & 0. 
\eeqn
Equality holds if and only if  $r_J({\cal F}) \leq 2$. And when
equality holds then $depth~G({\cal F}) \geq  3 - 2 + 1 = 2$
\cite[Theorem~3.7]{marley}.  

Similarly, 
\beqn
      \sum_{i \geq 1}h_i (J,{\cal F})
&=&    h_1 (J,{\cal F} ) - h_2 (J,{\cal F} ) +  h_3 (J,{\cal F} )\\
&=&   \sum_{n \geq 2}  
      \left[
      \ell \left( H_1 (C.( n)) \right) 
-     \ell \left( H_2 (C.(n)) \right) 
+      \ell \left( H_3 (C.( n)) \right)  \right]\\
&=&   \sum_{n \geq 2}  \left[
      \ell \left( H_0 (C.( n)) \right)
-     \ell \left( \f{R}{{\olin {I^{n}}}} \right)
+   3 \ell \left( \f{R}{ {\olin {I^{n-1}}}} \right)
-   3 \ell \left( \f{R}{{\olin {I^{n-2}}}} \right)
+     \ell \left( \f{R}{ {\olin {I^{n-3}}}} \right) \right]\\
&=&   \sum_{n \geq 2} \ell \left( \f{R}{J + {\olin {I^{n}}}} \right)
 -    \sum_{n \geq 3}  e(J)
-     \ell \left( \f{R}{{\olin {I^{2}}}} \right)
+    3\ell \left( \f{R}{{\olin I}} \right) \\
&=&   -\sum_{n \geq 2} \ell \left( \f{J + {\olin {I^{n}}}}{J} \right)
+     \sum_{S \succeq R}{o (I^S) \choose 3} [S/M_S : R/M]\\
&\geq& 0
\eeqn
If the above inequality is an equality, then  $depth~G({\cal F})
      \geq  3$. Conversely,  if $depth~G({\cal F}) \geq  3$ then there
      exists a minimal reduction $J$ such that
      above inequality is an
      equality. But if
      $depth~G({\cal F}) 
      \geq  3$, then $J \olin{I ^n} = \olin{I ^{n+1}}$ for all $n
      \geq 2$, i.e., $J + \olin{I^n} = J$ for all $n \geq 3$. Hence
\beqn
 \sum_{S \succeq R}{o (I^S) \choose 3} [S/M_S : R/M]
= \ell \left( \f{J + {\olin {I^{2}}}}{J} \right).  
\eeqn

Suppose $o(I) \leq 2$, then  $o( (I^{S})) \leq 2$ for all $S \succeq
      R$. Hence

\beqn
      \sum_{i \geq 1}h_i (J,{\cal F})
= - \ell \left( \f{J + {\olin {I^{2}}}}{J} \right)
\geq 0.
\eeqn
Since the length of the module appearing  above is
      positive, it is equal to zero. Hence  $\sum_{i \geq 1}h_i
      (J,{\cal F})=0$ which implies that  $depth~G({\cal F}) \geq  3$,
      i.e. $G({\cal F})$ is Cohen-Macaulay. 
This proves (3). 

If $r_J({\cal F}) = 1$ for some minimal reduction of $J$ of $I$, then
$G({\cal F})$ is Cohen-Macaulay.  Hence we have
\beqn
    0
=  \ell \left( \f{ J +  \olin{I^2}}{J} \right)
=  \sum_{S \succeq R}{o(I_S) \choose 3} [S/M_S : R/M]
\eeqn
This implies that that $o(I) \leq 2$. 
\qed

\section{Examples}

We end this paper with a few examples which will clarify our results and 
the assumptions we use. 

If $\mu(I) = d-1 + o(I)^{d-1}$, then $F(I)$
is Cohen-Macaulay if and only if there exists an ideal $J \subseteq I$
generated by $d$ elements such that $JI = I^2$ \cite{c-r-v}. 
Here we demonstrate an
example of a monomial ideal $I$ whose fiber cone is not
Cohen-Macaulay. This example also shows that $MI$ is not integrally
closed. It is easy to see that $\mu(I) = 11 > 10$.

\bex
\label{ex-7.1}
Let $I = (x^4, x^3y, x^2z, x^2y^2, xy^2z, xyz^2, xz^3, y^3, y^2z^2, yz^3,
z^5)$ be an ideal in the polynomial ring $k[x,y,z]$. Assume that 
$k$
is an algebraically closed field.  Then 
\been
\item
$I$ is a finitely supported ideal. 

\item
$J = (x^4 + yz^3, x^2z, y^3 + z^5 )$ is a minimal
reduction of $I$ and $r_{J}(I) = 2$. 

\item
Since $\mu(I) = 11 = o(I)^2 + 2$.  Since $r({\cal F}) = r(I) =
2$, $F(I)$ is not Cohen-Macaulay \cite[Corollary~2.5]{c-r-v}.  The
Hilbert  series of $F(I)$ is 
$$
H(F(I), t) = \f{1+8t}{(1-t)^3}. 
$$
Notice that $MI \not = \olin{MI}$ since $xy^2z \in  \olin{MI} \sms
MI$. 
\item
$I^n = \olin{I^n}$ for all $n \geq 1$ and $J \cap I^n = I^{n+1}$ for
all $n \geq 2$. 

\item
${\cal B} (I) 
=  \{ o(I^S) : S \succeq R \} 
=  \{ 3,2,1,1,1,1,1\}$. 
The Hilbert function 
$H_I(n) = \ell (R/ I^n)$ 
is equal to the Hilbert polynomial $P_I(n)$ for all $n \geq 0$.
In particular, 
\beqn
    \ell \left( \f{R}{I^n} \right)
&=& \ell \left( \f{R}{\olin{I^n}} \right) \\
&=&   {3n + 2 \choose 3}
+   5 {n  + 2 \choose 3}
+     {2n + 2 \choose 3}\\
&=& 40{n  + 2 \choose 3}
-   22{n  + 1 \choose 2}
+     {n      \choose 1}.
\eeqn

\item
By \cite[Theorem~17]{itoh}, $G({\cal F}) = G(I)$ is Cohen-Macaulay.

\item
Since 
$
  \ell(J + I^2/ J) 
= \sum_{S\succeq R} {o (I^S) \choose 3} 
= 1
$, 
by  Theorem~\ref{dim3},  $G({\cal F}) = G(I)$ is Cohen-Macaulay. 
\eeen
\qed
\eex

We now show that Theorem~\ref{main1} does not hold true if $I$ is 
not a monomial ideal. 

\bex
\label{extwo}
{\em 
(The author  is  very grateful to Oliver Piltant for bringing this
example into light.)
Consider the following example. Let $k$  be an algebraically closed
field. Let $R = k[x,y,z]$ and  let 
$I = (z^3, y^3 - x^2 z, y^2 z^2, x y z^2, 
      x^2 z^2, x y^2 z, x^2 y z, x^3 z, x^2 y^2, x^3 y, x^4)
$ 
be an ideal of $R$. Then it is easy to verify that 

\been
\item
$I$ is complete and $\ell (R/ I) = 18$;

\item
The strict transform of $I$ in $S = R[y/x, z/x]$ is 
$I^S = (x, (y/x)^3 - (z/x), (z/x)^3)$ and 
$\ell (S/ I^S) =9$. 
But 
$
\displaystyle { {o(I) + 2 \choose 3} = 10}$. 
Hence
$
\ell ( R^1 \sigma_{*} ( I {\cal O}_{X({\cal C}_I)})) = 1
$. 

\item
Note that if  $J = (z^3, y^3 - x^2 z,  x^4)$, then $J$ is generated by
a system of parameters in $I$ and $JI =I^2$. Hence    $F(I)$ is
Cohen-Macaulay \cite[Corollary~2.5]{c-r-v}.  Also $G(I)$ is
Cohen-Macaulay. 
\eeen
\qed}
\eex

We now present an infinite class of ideals where  $\ell ( R^i
\sigma_{*} ( I {\cal O}_{X({\cal C}_I)})) >0$. This is a
generalization of the example in \cite{huc-huneke}. 

\bex
\label{ex-7.2}
{\em Let $R = \comx[x,y,z]$ where $x,y,z$ are variables and let 
$$
I = (x^{r+1}, (x,y,z)(y^r + z^r), (x,y,z)^{r+2}).
$$
Then $I$ is a finitely supported complete ideal. Put 
$S = R[x/y, z/y]$. Then 
\beqn
     \ell \left( \f{R}{I} \right)
&=&   {r + 3 \choose 3} 
+     {r + 3 \choose 2} -4 \\
     \ell \left( \f{S}{I^S} \right)
+     {o(I) + 2 \choose 3 }
&=&   2 {r + 1 \choose 2} 
+       {r + 3 \choose 3} \\
     \ell ( R^i
     \sigma_{*} ( I {\cal O}_{X({\cal C}_I)}))
&=&  \ell \left( \f{R}{I} \right)
-    \ell \left( \f{S}{I^S} \right)\\
&=&     {r - 1 \choose 2}.
\eeqn
Using the argument on the lines in \cite{huc-huneke} one can show
that $G(I)$ is not Cohen-Macaulay for all $r \geq 3$. In particular, 
$x^3(y^r+z^r)^3 \in J \cap I^3 \sms JI^2$ where
$J = (x^{r+1}, z(y^r + z^r), y(y^r + z^r) + yz^{r+1})$ is  a minimal
reduction of $I$.
\qed}
\eex

We end this paper with the following example. 
\bex
\label{order}
{\em
Let $k$ be an algebraically closed field. Let $R=k[x,y,z]$. Let $M= (x,y,z)$,
$I_1 = (x, y^2, yz, z^2)$ and $I_2 = (x^2, y, z)$ and put $I=I_1 I_2$.
Then $I$ is not a finitely supported ideal. 
We demonstrate the fact that Corollary~\ref{mixed-mult}(2) and 
Theorem~\ref{generators}(2),(3)
does not hold true if the ideal is not finitely supported.

By \cite{teissier} and \cite{rees}
\beqn
e_2(M|I_1) &=& e(x,y,z^2) = 2 \\
e_2(M|I_2) &=& e(x,y,z) = 1 \\
e(M|I_1|I_2) &=&  e_{1,1,1}(M,I_1,I_2) = e(x,y,z) = 1
\eeqn
Hence  
\beqn
e_2(M|I_1 I_2) = e_2(M|I_1) + e_2(M|I_2) + 2 e(M|I_1|I_2)= 5 > 4 = o(I) \\
\mu(I) = 7 > 6 = o(I)^2 + (3-1). 
\eeqn
It is also easy to verify that $MI_1$ is integrally closed, but 
$$
\mu(I_1) = 4 > 3 = {o(I_1) + 2 \choose 2}. 
$$
}
 \eex
\qed

\end{document}